\documentclass[12pt]{amsproc}

\usepackage{amscd}

\newtheorem{lemma}{Lemma}

\newtheorem{theorem}[lemma]{Theorem}
\newtheorem{prop}[lemma]{Proposition}

\def\proof{\par\medskip\noindent {\sc Proof. }}
\def\remark{\par\medskip \noindent {\sc Remark. }}
\newcommand{\proofof}[1]{\par\medskip\noindent {{\sc Proof of #1. }}}
\def\qed{\hfill $\square$\bigskip}

\newfont{\bbf}{msbm10 at 12pt}

\def\N{\mathbb N}
\def\C{\mathbb C}
\def\Z{\mathbb Z}

\newcommand{\capt}[1]{{\bf(#1)}\\}

\renewcommand{\Re}{{\rm Re}}

\newcommand{\s}{{\underline s}}

\newcommand{\diam}{{\mbox{\rm diam}}}
\newfont{\Euler}{eusm10 at 12pt}
\newcommand{\Sym}{ \mathcal{S} }
\newcommand{\sm}{\setminus}
\newcommand{\hide}[1]{}

\newcommand{\Ray}{G_\s}
\newcommand{\wt}{\widetilde}

\newcommand{\rayK}{g_{\kappa,\s}}
\newcommand{\rayKt}{g_{\kappa(t),\s}}
\newcommand{\raysigK}{g_{\kappa,\sigma(\s)}}
\newcommand{\raysign}{g_{\kappa,\sigma^n(\s)}}
\newcommand{\raysignt}{g_{\kappa(t),\sigma^n(\s)}}
\newcommand{\raysigm}{g_{\kappa,\sigma^m(\s)}}
\newcommand{\raysigmn}{g_{\kappa,\sigma^{m+n}(\s)}}

\title[Hausdorff Dimension of Parameter Rays]
{Hausdorff Dimension of Exponential Parameter Rays and Their Endpoints}
\author{Mihai Bailesteanu}
\author{Horia Vlad Balan}
\author{Dierk Schleicher}

\address{Department of Mathematics, Malott Hall, Cornell University, Ithaca, NY 14853, USA}
\email{mb452@cornell.edu}

\address{Ming Hsieh Department of Electrical Engineering, University of Southern California,
Hughes Aircraft Electrical Engineering Center, Los Angeles, CA 90089-2560}
\email{vlad.gm@gmail.com}

\address{School of Engineering and Science, Jacobs University Bremen (formerly International University Bremen), Postfach 750 561, D-28725 Bremen, Germany}
\email{dierk@jacobs-university.de}

\begin{document}



\begin{abstract}
We investigate the set $I$ of parameters $\kappa$ for which the singular value of $z\mapsto e^z+\kappa$ converges to $\infty$. The set $I$ consists of uncountably many parameter rays, plus landing points of some of these rays \cite{FRS}. We show that the parameter rays have Hausdorff dimension $1$, which implies \cite{QiuDimTwo} that the ray endpoints in $I$ alone have dimension $2$. Analogous results were known for dynamical planes of exponential maps \cite{Karpinska,SZ1}; our result shows that this also holds in parameter space.
\end{abstract}

\maketitle


\section{Introduction}
\label{sec:intro}

We study the space of exponential maps $E_\kappa\colon z\mapsto e^z+\kappa$: each map $E_\kappa$ has a unique singular value, which is $\kappa$ (equivalently, one often uses the parametrization $z\mapsto\lambda e^z$ with the unique singular value $0$). This space of exponential maps is decomposed into the locus of \emph{structurally stable maps}, i.e., those maps $E_\kappa$ which are topologically conjugate to $E_{\kappa'}$ for all $\kappa'$ sufficiently close to $\kappa$, and the complementary \emph{bifurcation locus} $B$. Many structurally stable exponential maps are \emph{hyperbolic}, i.e., the orbit of the singular value $\kappa$ converges to an attracting periodic orbit. Hyperbolic exponential maps are completely classified \cite{ExpoAttr}; a fundamental conjecture states that all structurally stable maps are hyperbolic.

The bifurcation locus $B$ is very complicated. It contains the \emph{escape locus}, which is the set 
\[
I:=\left\{
\begin{minipage}[c]{68mm}{$\kappa\in\C\colon$ 
the orbit of the singular value converges to $\infty$
under iteration of $E_\kappa$}
\end{minipage}
\right\}\,\,.
\]
The set $I$ is completely classified; see below. It contains uncountably many curves called \emph{parameter rays}. It is conjectured that every $\kappa\in B$ is either on a parameter ray or the unique limit point of finitely many parameter rays, so a complete description of the bifurcation locus could be given in terms of $I$ (this would imply the previous conjecture \cite[Corollary~A.5]{RS}).

The situation is analogous for the space of quadratic polynomials: the bifurcation locus is the boundary of the famous Mandelbrot set, all hyperbolic components are explicitly classified, and the fundamental conjectures are that all structurally stable components are hyperbolic and that every point in the bifurcation set is the unique limit point of a positive finite number of parameter rays (the second conjecture is equivalent to \emph{local connectivity of the Mandelbrot set} \cite{ExpoBif}, and it implies the first \cite{Orsay,Fibers2}). 

The set $I$ was studied in detail in \cite{FRS}, building on earlier work in \cite{dghnew2,FS}. A parameter ray is a maximal injective curve $\Ray\colon (t_\s,\infty)\to I$ with $\Ray(t)\to\infty$ as $t\to\infty$ for some $t_\s\in[0,\infty)$. It was shown in \cite[Theorem~1.1]{FRS} that every path component of $I$ is either a unique parameter ray, or a unique parameter ray $\Ray$ together with its unique endpoint $\kappa_\s=\lim_{t\searrow 0}\Ray(t)$ (in \cite{FS,FRS} a preferred parametrization is described that determines $t_\s$ uniquely).

Let $I_R$ be the union of all parameter rays, and $I_E:=I\sm I_R$ be the set of endpoints in $I$. Our main result is the following.

\begin{theorem}\capt{Hausdorff Dimension of Parameter Rays and Endpoints}
\label{Thm:ParaRaysEndpoints}%
The Hausdorff dimension of the sets of parameter rays $I_R$ and their escaping endpoints $I_E$ satisfies $\dim_H(I_R)=1$ and $\dim_H(I_E)=2$. 
\end{theorem}

Results like this are sometimes called a ``dimension paradox'': every path component of $I$ contains a curve in $I_R$ and at most a single endpoint in $I_E$, so one might think that ``most'' points in $I$ are in $I_R$, i.e., on a ray. Topologically, this is correct, but nonetheless the entire Hausdorff dimension sits in $I_E$. This is not the first time such a phenomenon is observed in transcendental dynamics: in the dynamic plane of every exponential map $E_\kappa$, the set of escaping points $I_\kappa:=\{z\in\C\colon E_\kappa^{\circ n}(z)\to\infty \mbox{ as }n\to\infty\}$ is known to have a very similar structure \cite{Karpinska,SZ1}, and stronger results hold for maps of the form $z\mapsto \pi\sin z$: here, the set of endpoints can have full planar Lebesgue measure and more \cite{CosDuke,Monthly}. Our result is the first which establishes this ``dimension paradox'' in parameter space. Qiu~\cite{QiuDimTwo} proved that $I$ has Hausdorff dimension $2$ and zero planar Lebesgue measure, so for us it suffices to focus on the set $I_R$ of parameter rays and show that it has Hausdorff dimension $1$; the claim $\dim_H(I_E)=2$ is then immediate.
(Note that for the Mandelbrot set, it is known that the boundary has Hausdorff dimension $2$ \cite{MitsuHausdorff}; in this case, the union of the parameter rays forms an open subset of $\C$, so it has dimension $2$.)

The results in this paper are the main results of the Bachelor theses \cite{VladThesis} and \cite{MihaiThesis}. They owe a lot to discussions with Markus F\"orster; his contributions are gratefully acknowledged. We would also like to thank Lasse Rempe and two anonymous referees for helpful comments.

{\bf Notation.}
We denote the $n$-th iterate of $E_\kappa$ by $E_\kappa^{\circ n}$ and write $E^n(\kappa):=E_\kappa^{\circ n}(\kappa)$ for integers $n\ge 0$. Then every $E^n$ is a transcendental entire function and $I:=\{\kappa\in\C\colon E^n(\kappa)\to\infty$ as $n\to\infty\}$.

\section{Standard Squares and Conformal Mappings}
\label{Sec:Conformal}

{The fundamental idea of our proof will be similar to Hausdorff dimension estimates for dynamic rays in the dynamic planes of exponential maps as in \cite{Karpinska,SZ1}; however, the maps $E^n$ are not iterates, so the main task is to establish good local mapping properties of these maps, close to what one has for iterates. This will be done in this section.}

If $\Lambda\subset\C$ is a domain so that $E^n\colon\Lambda\to E^n(\Lambda)=:V$ is a conformal isomorphism, then for every $k\ge 0$ this defines a holomorphic map $E^{n,n+k}=E^{n+k}\circ (E^n)^{-1}\colon V\to\C$.

\begin{lemma}\capt{Univalent Extra Iterate}
\label{Lem:UnivalentExtraIterate}%
Suppose $\Lambda,V\subset\C$ are such that $E^n\colon\Lambda\to V$ is a conformal isomorphism, $\Re(V)>\xi>1$ and $|(E^n)'|>2$ for all $\kappa\in\Lambda$; moreover, suppose that $V$ is convex and contained in a horizontal strip of height $\pi/2$. Then $E^{n+1}\colon\Lambda\to E^{n+1}(\Lambda)$ is a conformal isomorphism with 
\[
|(E^{n+1})'(\kappa)|>e^\xi|(E^n)'(\kappa)|-1>2|(E^n)'(\kappa)| 
\]
for all $\kappa\in\Lambda$, 
and $E^{n,n+1}\colon V\to E^{n+1}(\Lambda)$ is a conformal isomorphism
with \[
|(E^{n,n+1})'(z)| > e^\xi-1
\]
for all $z\in V$.
\end{lemma}

\proof
In order to show that $E^{n+1}$ restricted to $\Lambda$ is a conformal isomorphism onto its image, all we need to check is injectivity of $E^{n+1}$ on $\Lambda$, or equivalently of $E^{n,n+1}$ on $V$.
We can write $E^{n,n+1}(z)=e^z+\kappa$ with $\kappa=(E^n)^{-1}(z)$. 

Suppose there are $\kappa_1,\kappa_2\in\Lambda$ with $E^{n+1}(\kappa_1)= E^{n+1}(\kappa_2)$. Set $z_j=E^n(\kappa_j)$ for $j=1,2$. Then $e^{z_1}+\kappa_1=e^{z_2}+\kappa_2$ or
\[
\kappa_2-\kappa_1 = e^{z_1}-e^{z_2} \,\,.
\]

Let $\gamma$ be the straight line segment connecting $z_1$ to $z_2$; we have $\gamma\subset V$ by convexity; since $(E^n)'>2$ on $\Lambda$, we have $|\kappa_2-\kappa_1|\le |z_1-z_2|/2$. Similarly, $|\exp'_V|>1$, so $\exp(\gamma)$ is a curve in $\C$ connecting $e^{z_1}$ and $e^{z_2}$ with length at least $|z_1-z_2|$. Since imaginary parts of $z\in V$ cannot change by more than $\pi/2$, the unit tangent vector of $\exp(\gamma)$ is always within a sector of width $\pi/2$ (or $90^\circ$) and thus $|e^{z_1}-e^{z_2}|\ge |z_1-z_2|/\sqrt 2$. As a result,
\[
|\kappa_2-\kappa_1| \le \frac 1 2 |z_1-z_2| \le \frac {\sqrt 2} 2 |e^{z_1}-e^{z_2}|
\]
and this is possible only if $\kappa_2=\kappa_1$.

Now we discuss the derivative of $E^{n+1}$. For this, 
\[
(E^{n+1})'=\exp\left(E^n\right) \cdot (E^n)' + 1
\]
implies
\[
|(E^{n+1})'|\ge \exp\left(\Re(E^n)\right) \cdot |(E^n)'| - 1
\ge e^\xi|(E^n)'|-1 > 2|(E^n)'|
\]
as claimed. Similarly,
\[
(E^{n,n+1})'(z) = e^z+d\kappa/dz = e^z+1/(E^n)'(\kappa)
\]
and so $|(E^{n,n+1})'(z)|>e^\xi - 1$.
\qed

We define a {\em standard square} $Q$ to be an open square of side length $\pi/2$ with sides parallel to the real and imaginary axes. The {\em double square} of a standard square is a square $\tilde Q\supset Q$ of side length $\pi$ with parallel sides and common center. We write $D_r(z)$ for the open disk of radius $r$ around $z\in\C$.

For real numbers $p>1$ and $\xi>0$, we define truncated parabola domains
\[ 
P_{p,\xi}:=\left\{z=x+iy\in\C\colon x>\xi \mbox{ and } |y|<x^{1/p} \right\}
\,\,.
\] 

\begin{lemma}\capt{Covering by Disks}
\label{Lem:CoveringDisks}%
Fix $p>1$, $\xi>0$ and $\Lambda\subset\C$ open. If  $\kappa\in I\cap\Lambda$ is such that $E^n(\kappa)\in P_{p,0}$ for all but finitely many $n$, and $|(E^n)'(\kappa)|\to\infty$ as $n\to\infty$, then there are an $N\in\N$, a neighborhood $U\subset\Lambda$ of $\kappa$ and a standard square $Q\subset P_{p,\xi}$ with center $E^N(\kappa)$ and double square $\wt Q$ so that $E^N\colon U\to \wt Q$ is a conformal isomorphism
and $E^n(\kappa)\in P_{p,\xi}$ for all $n\ge N$, and so that $|E^N(\kappa')|>2$ for all $\kappa'\in U$.
\end{lemma}
\proof
We may as well suppose that $e^\xi>33$.
Let $N_0\in\N$ be such that all $n\ge N_0$ satisfy
$E^n(\kappa)\in P_{p,\xi}$ and $|(E^n)'(\kappa)|>2$. Then 
$\kappa$ has a neighborhood $U_0$ so that $E^{N_0}\colon U_0\to E^{N_0}(U_0)$ is a conformal isomorphism. 

Let $r_0>0$ be the largest radius so that $D_{r_0}(E^{N_0}(\kappa))\subset E^{N_0}(U_0)$. If $r_0<\sqrt 2\pi/2$ then restrict $r_0$ if necessary so that $r_0\le \pi/4$. By Lemma~\ref{Lem:UnivalentExtraIterate}, the maps $E^{N_0+1}\colon U_0\to E^{N_0+1}(U_0)$ and $E^{N_0,N_0+1}(E^{N_0}(U_0))\to E^{N_0+1}(U_0)$ are conformal isomorphisms and we have $|(E^{N_0,N_0+1})'(E^N_0(\kappa))| > e^\xi-1$. By the Koebe $1/4$-theorem, there is a neighborhood $U_1\subset U_0$ of $\kappa_0$ so that $E^{N_0+1}(U_1)$ is a disk of radius $r_1\ge (e^\xi-1)r_0/4>8r_0$. Repeating this argument finitely many times, we obtain an index $N$ and a domain $U\ni\kappa$ so that $E^N\colon U\to \wt Q$ is a conformal isomorphism, where $\wt Q$ is a double square. The condition $|(E^N)'|>2$ on $U$ can be assured by finitely many extra iterations if necessary, because derivatives grow uniformly by Lemma~\ref{Lem:UnivalentExtraIterate}.
\qed

\section{Hausdorff Dimension Estimates}
\label{Sec:Dimension}

For bounded open sets $\Lambda\subset\C$, we are interested in the set 
\[
I_{p,\Lambda}:=\left\{
\begin{array}{l}
\kappa\in\Lambda\cap I\colon  \mbox{$|(E^n(\kappa))'|\to\infty$ as $n\to\infty$} \\
\mbox{and $E^n(\kappa)\in P_{p,0}$ for all sufficiently large $n$}
\end{array}
\right\}
\,\,.
\]
Set also 
\[
I_{p,\xi,\Lambda}^N:=\{\kappa\in I_{p,\Lambda}\colon
E^n(\kappa)\in P_{p,\xi} \mbox{ for all $n\ge N$}\}
\,\,.
\]

\begin{prop}\capt{Hausdorff Dimension Estimate}
\label{Prop:DimEstimate}%
Fix $p >1$ and an integer $N\ge 0$. Suppose $Q\subset\C$ is a standard square with double square $\wt Q$ and $\wt\Lambda\subset\C$ is such that $E^N\colon \wt\Lambda\to\wt Q$ is a conformal isomorphism. Suppose also that $|(E^N)'|>2$ on $\wt\Lambda$. Set $\Lambda:=(E^N)^{-1}(Q)\cap\wt\Lambda$, and let $M>0$ be such that $\Lambda\subset D_M(0)$.
Then $\dim_H(I_{p,\xi,\Lambda}^N) \le 1+1/p$ provided $\xi$ is sufficiently large depending only on $p$ and $M$.
\end{prop}
\proof
Let $\xi_0$ be such that $\Re(Q)\subset(\xi_0,\xi_0+\pi/2)$. We may suppose that $\xi_0\ge \xi-\pi/2$ because otherwise $I_{p,\xi,\Lambda}^N=\emptyset$.

By Lemma~\ref{Lem:UnivalentExtraIterate}, $E^{N,N+1}\colon Q\to E^{N,N+1}(Q)=:W$ is a conformal isomorphism; we can write $E^{N,N+1}(z)=\exp(z)+\kappa$ with $\kappa\in\Lambda$. The set $\exp(Q)$ is contained in an annulus between radii $e^{\xi_0}$ and $e^{\pi/2} e^{\xi_0}$, and $E^{N,N+1}(Q)\cap P_{p,\xi}$ has real parts between $\xi_1:=e^{\xi_0}/2$ and $e^{\pi/2} e^{\xi_0}+M$ (provided $\xi$ is sufficiently large). Consequently the imaginary parts in $E^{N,N+1}(Q)\cap P_{p,\xi}$ have absolute values at most $e^{\pi/2p}e^{\xi_0/p}+1$ (again for sufficiently large $\xi$). Therefore, $W\cap P_{p,\xi}$ can be covered by at most 
\[
N(\xi_0):=\left(e^{\pi/2} e^{\xi_0}+M-e^{\xi_0}/2\right) \cdot 2(e^{\pi/2p}e^{\xi_0/p}+1) (\pi/2)^{-2}\le C e^{\xi_0(1+1/p)}
\] 
standard squares for some universal constant $C>0$; denote these $N(\xi_0)$ standard squares by $Q_{1,i}$ for $i=1,2,\dots,N(\xi_0)$.

Denoting the double squares of $Q_{1,i}$ by $\wt Q_{1,i}$, we have $E^{N,N+1}(\wt Q)\supset \bigcup \wt Q_{1,i}$ (this just needs the fact that $|(E^{N,N+1})'|$ is large on $\wt Q$). We can thus pull back  the $Q_{1,i}$ under $E^{N,N+1}$ and obtain a covering of $Q\cap (E^{N,N+1})^{-1}(P_{p,\xi})$ with $N(\xi_0)$ open sets $U_{1,i}=W_{1,i}$ so that each $U_{1,i}$ has a neighborhood $\wt U_{1,i}$ for which the restriction $E^{N,N+1}\colon \wt U_{1,i}\to \wt Q_{1,i}$ is a conformal isomorphism. By the Koebe distortion theorem, the restrictions $E^{N,N+1}\colon U_{1,i}\to Q_{1,i}$ have uniformly bounded distortions, and their derivatives are at least $e^{\xi_0}-1$. 
Note that for any $d$
\[
\sum_i(\diam\,U_{1,i})^d \le C' N(\xi_0)(\pi / (e^{\xi_0}-1))^d
\le C'' e^{\xi_0(1+1/p-d)}
\,\,,
\]
where $C'$ and $C''$ are universal constants; $C'$ measures the distortion of $E^{N,N+1}\colon U_{1,i}\to Q_{1,i}$. In particular, if $d>1+1/p$ is fixed and $\xi$ is sufficiently large, then 
\begin{equation}
\sum_i(\diam\,U_{1,i})^d\le (\diam\,Q )^d \,\,.
\label{Eq:HausdorffMeasureDecreases}
\end{equation}

This argument can be repeated: each standard square $Q_{1,i}$ has real parts at least $\xi_1=e^{\xi_0}/2\gg\xi_0\ge\xi-\pi/2$, so $Q_{1,i}\cap (E^{N+1,N+2})^{-1}(P_{p,\xi})$ can be covered by at most $N(\xi_1)$ open sets $W_{2,i'}$ so that the image sets $E^{N+1,N+2}(W_{2,i'})$ are in turn standard squares $Q_{2,i'}$ at real parts at least $\xi_2\gg\xi_1$, and so on. The sets $W_{2,i'}$ can be pulled back under $E^{N,N+1}$ and yield a covering $U_{2,i'}$ of the set
\[
\left\{ z\in Q\colon E^{N,N+1}(z)\in P_{p,\xi} \mbox{ and } E^{N,N+2}(z)\in P_{p,\xi}\right\} \,\,.
\]
Inductively, we obtain a family of coverings $U_{n,i}$ for every $n\ge 1$, and the set they cover is
\[
\hat Q:=\left\{ z\in Q\colon E^{N,N+n}(z)\in P_{p,\xi} \text{ for all $n\ge 1$} \right\} \,\,.
\]
Each $U_{n,i}\subset\wt Q$ is such that $E^{N,N+n}\colon U_{N,i}\to Q_{n,i}$ is a conformal isomorphism, where $Q_{n,i}$ are standard squares at real parts at least $\xi_n$, and $E^{N,N+n}$ extends to a conformal isomorphism $\wt U_{n,i}\to\wt Q_{n,i}$ where $\wt Q_{n,i}$ is the double square of $Q_{n,i}$. The sets $W_{n,i}:=E^{N,N+n-1}(U_{n,i})$ cover those points $z$ in the standard squares $Q_{n-1,i'}$ of the previous generation for which $E^{N+n-1,N+n}(z)\in P_{p,\xi}$.

We show that $\dim_H(\hat Q)\le 1+1/p$. Indeed, for $d>1+1/p$ and sufficiently large $\xi$ it follows that the maps $E^{N,n}\colon U_{n,i}\to Q_{n,i}$ are conformal isomorphisms with uniformly bounded distortions, so (\ref{Eq:HausdorffMeasureDecreases}) becomes
\begin{equation}
\sum_{i' }(\diam\, U_{n+1,i'})^d < \sum_i (\diam\, U_{n,i})^d 
\,\,.
\label{Eq:HausdorffMeasureDecreases2}
\end{equation}
However, for fixed $d>1+1/p$, we cannot be sure that $\xi$ is sufficiently large. But the $\xi_n$ grow exponentially fast, and (\ref{Eq:HausdorffMeasureDecreases}) holds for all sufficiently large $n$.

Since $E^{N,N+n}\colon U_{n,i}\to Q_{n,i}$ are conformal isomorphisms with boun\-ded distortions and derivatives tending to $\infty$ as $n\to\infty$, it follows that $\sup_i\diam\, U_{n,i}\to 0$ as $n\to\infty$. The family of covers $U_{n,i}$ proves that $\dim_H(\hat Q)\le 1+1/d$. 
Finally, $E^N\colon \Lambda\to Q$ is a conformal isomorphism with $E^N(I_{p,\xi,\Lambda}^N)\subset \hat Q$. Therefore, $\dim_H(I_{p,\xi,\Lambda}^N)\le 1+1/d$ as well.
\qed

\begin{theorem}\capt{Hausdorff Dimension of Parameter Rays}
\label{Thm:HausdorffDimParaRays}%
For every $p>1$ and every bounded open $\Lambda\subset\C$, we have $\dim_H(I_{p,\Lambda})\le 1+1/p$.
\end{theorem}
\proof
Choose $\xi>0$ depending on $\Lambda$ and $p$ as in Proposition~\ref{Prop:DimEstimate}. Pick some $\kappa\in I_{p,\Lambda}$. By Lemma~\ref{Lem:CoveringDisks}, there are an $N\in\N$, a neighborhood $U\subset\Lambda$ of $\kappa$ and a standard square $Q\subset P_{p,\xi}$ with double square $\wt Q$ so that $E^N\colon U\to\wt Q$ is a conformal isomorphism with $|(E^n)'|>2$, $E^N(\kappa)$ is the center of $Q$ and $E^n(\kappa)\in P_{p,\xi}$ for all $n\ge N$. Then $\kappa\in I_{p,\xi,U}^N$ and by Proposition~\ref{Prop:DimEstimate}, $\dim_H(I_{p,\xi,U}^N)\le 1+1/p$. 

Since $\Lambda$ has countable topology and $N$ is from a countable set, 
$I_{p,\Lambda}$ is contained in the countable union of sets of dimension at most $1+1/p$, and the claim follows.
\qed

In order to prove Theorem~\ref{Thm:ParaRaysEndpoints}, we need to introduce \emph{parameter rays} $\Ray$ as introduced in \cite{FS,FRS} together with a particular parametrization. For every $\s\in\Sym:=\Z^\N$, there is a well-defined \emph{minimal potential} $t_\s\in[0,\infty]$ and an injective curve $\Ray\colon(t_\s,\infty)\to\C$ with the following properties:
\begin{itemize}
\item
for every $\kappa\in I_R$ there is a unique external address $\s\in\Sym$ and a unique potential $t>t_\s$ with $\kappa=\Ray(t)$;
\item
if $\kappa=\Ray(t)$, then for $E_\kappa$ the asymptotics of the singular orbit can be expressed using $F(t):=e^t-t$ and $\s=s_1s_2s_3\dots\in\Sym$:
\begin{equation}
E^n(\kappa)=E_{\kappa}^{\circ n}(\kappa)=F^{\circ n}(t)+2\pi i s_{n+1}+O(1) \qquad \mbox{as $n\to\infty$}
\,\,.
\label{Eq:ParaRayIterate2}
\end{equation}
\end{itemize}
More precisely, \cite[Corollary~3.2 and Proposition~2.2]{FRS} justify the partition $I=I_R\cup I_E$ of escaping parameters into parameter rays and escaping endpoints, so that $I_R$ is exactly the set of parameter rays as classified in \cite[Theorem~3.7]{FS}. By \cite[Theorem~3.12]{FS}, every parameter ray is a $C^1$-curve $\Ray\colon(t_\s,\infty)\to\C$ for a well-defined $\s\in\Sym$ and $t_\s\ge 0$, so that all parameter rays are injective and disjoint curves with $\Ray'(t)\neq 0$ for all $t>t_\s$. Not all $\s\in\Sym$ actually occur as external addresses; for the others we have $t_\s=\infty$. (Note that the maps in \cite{FS} are parametrized as $z\mapsto \exp(z+\kappa)$, rather than our $z\mapsto \exp(z)+\kappa$; but these maps are conjugate by translation.)

\begin{lemma}\capt{Derivative of Parameter Rays}
Parameter rays satisfy $(d/d\kappa)E^n(\kappa)\to\infty$ and $(d/dt)E^n(\Ray(t))\to\infty$  as $n\to\infty$.
\label{Lem:ParaRaysDiffEstimate}
\end{lemma}

\remark
This result is a key ingredient in the proof \cite{ExpoBif,ExpoPara} that exponential parameter space is not locally connected at any point on any parameter ray.
\proof
We will prove the result using rays in the dynamic planes: \emph{dynamic rays} for $E_\kappa$ are curves in $\C$ consisting entirely of \emph{escaping points}, i.e., points $z$ with $E_\kappa^{\circ n}(z)\to\infty$ as $n\to\infty$.

The existence of dynamic rays was shown in \cite[Theorem~4.2]{SZ1}: 
For every $\kappa=\Ray(t)$, there are dynamic rays $\rayK\colon(t_\s,\infty)\to\C$ for every $\s\in\Sym$ and $t_\s\in[0,\infty]$ (these curves are empty if $t_\s=\infty$); however, there are exceptions if some dynamic ray contains the singular value (there are no preimages of the singular ray). In \cite[Theorem~3.7]{FS}, parameter rays are defined so that $\kappa=\Ray(t)$ if and only if $\rayK(t)$ equals the singular value $\kappa$: so we are exactly in the situation where exceptions to the existence of dynamic rays occur. However, it is shown in \cite[Section~3]{FS} that for $\kappa_0=\Ray(t)$ there is a neighborhood $\Lambda$ of $\kappa_0$ in parameter space so that $\rayK(t)$ is defined for all $\kappa\in\Lambda$. 

By \cite[Theorem~4.2]{SZ1}, dynamic rays satisfy 
$E_\kappa(\rayK(t))=\raysigK(F(t))$, where $\sigma$ is the left shift on the sequence $\s=s_1s_2s_3\dots$. By \cite[Proposition~4.5]{SZ1} we have the asymptotics
\begin{equation}
E_\kappa^{\circ n}(\rayK(t))=\raysign(F^{\circ n}(t))=F^{\circ n}(t)+2\pi i s_{n+1}+o(1)\to\infty
\label{Eq:OrbitAsymptotics}
\end{equation}
as $n\to\infty$. We thus consider
\begin{eqnarray*}
\frac{d}{dt} E^n(\Ray(t))
&=& \frac{d}{dt} E_{\kappa(t)}^{\circ n}(\rayKt(t))
= \frac{d}{dt} \raysignt(F^{\circ n}(t))
\\
&=&
\frac{\partial}{\partial t} \raysign(F^{\circ n}(t)) + \frac{\partial}{\partial \kappa} \raysign(F^{\circ n}(t)) \frac{d\kappa}{dt}
\\
&=& \raysign'(F^{\circ n}(t))\cdot \frac{d}{dt} F^{\circ n}(t) + \frac{\partial}{\partial \kappa} \raysign(F^{\circ n}(t)) \cdot \Ray'(t)
\,\,.
\end{eqnarray*}

We clearly have $dF^{\circ n} (t)/dt\to\infty$ as $n\to\infty$. We will now show that $\raysign'(F^{\circ n}(t))\to 1$ and $\frac{\partial}{\partial \kappa} \raysign(F^{\circ n}(t)) \to 0$; since $\Ray'(t)\in\C$ is a fixed number, this will imply $(d/dt)  E^n(\Ray(t))=\infty$ and thus prove the first claim of the lemma.

By \cite[Proposition~4.6]{FS}, we have for $t>t_\s$
\begin{equation}
\rayK'(t)=\prod_{m=1}^\infty \frac{F^{\circ m}(t)+1}{\raysigm(F^{\circ m}(t))}
\label{Eq:ProductDerivative}
\end{equation}
and thus
\[
(\raysign)'(F^{\circ n}(t))
=\prod_{m=1}^\infty \frac{F^{\circ (m+n)}(t)+1}{\raysigmn(F^{\circ (m+n)}(t))}
=\prod_{m=n+1}^\infty \frac{F^{\circ m}(t)+1}{\raysigm(F^{\circ m}(t))}
\,\,.
\]
The claim $(\raysign)'(F^{\circ n}(t))\to 1$ follows directly from convergence of (\ref{Eq:ProductDerivative}) (as shown in \cite{FS}).

For the last limit, we have to review how dynamic rays are defined: for $m\ge 0$, we define $g_{\kappa,\s}^0(t):=t$ and $g_{\kappa,\s}^{m+1}(t):=L_{s_1}(g_{\kappa,\sigma(\s)}^m(F(t)))$, where $L_s(z):=\log(z-\kappa)+2\pi i s$ is an inverse branch of $E_\kappa=e^z+\kappa$  (for $s\in\Z$). In \cite{SZ1}, dynamic rays are constructed as $\rayK(t)=\lim_{m\to\infty}g_{\kappa,\s}^m(t)$: in \cite[Proposition~3.4]{SZ1}, it is shown that the $g_{\kappa,\s}^m$ converge uniformly to a limiting curve $\rayK$ for sufficiently large $t$ (``on ray tails''), so that $\rayK(t)$ depends holomorphically on $\kappa$ for fixed $t$. Moreover, the convergence is locally uniform in $\kappa$ and all functions are holomorphic, so $(d/d\kappa)\rayK(t)=\lim_{m\to\infty}(d/d\kappa)g_{\kappa,\s}^m(t)$ (recall that in \cite{SZ1}, the parametrization $\exp(z+\kappa)$ with inverse $L_s(z)=\log(z)-\kappa+2\pi is$ is used; these are conjugate by translation). 

This construction is extended in \cite[Theorem~4.2]{SZ1} to entire dynamic rays, i.e., for all $t>t_\s$: a point $z$ is on a dynamic ray if $E_\kappa^{\circ n}(z)$ is on a ray tail for sufficiently large $n\ge 0$. Since we are interested in the limit $\frac{\partial}{\partial \kappa} \raysign(F^{\circ n}(t))$ as $n\to\infty$, we may restrict to sufficiently large $n$ so that we are always on ray tails. We need to prove
\[
0=\lim_{n\to\infty}\frac{\partial}{\partial \kappa} \raysign(F^{\circ n}(t)) =
\lim_{n\to\infty}\lim_{m\to\infty} \frac{\partial}{\partial\kappa} g_{\kappa,\sigma^n(\s)}^m(F^{\circ n}(t))  \,\,.
\]
From $g_{\kappa,\sigma^n(\s)}^{m+1}(F^{\circ n}(t))=L_{s_{n+1}}(g_{\kappa,\sigma^{n+1}(\s)}^m(F^{\circ(n+1)}(t)))$ we obtain the recursive relation
\[
\frac{\partial}{\partial\kappa} g_{\kappa,\sigma^n(\s)}^{m+1}(F^{\circ n}(t))=
\frac{\frac{\partial}{\partial\kappa}g_{\kappa,\sigma^{n+1}(\s)}^{m}(F^{\circ(n+1)}(t))-1}{g_{\kappa,\sigma^{n+1}(\s)}^{m}(F^{\circ(n+1)}(t))-\kappa} 
\]
starting with $(\partial/\partial\kappa)g_{\kappa,\sigma^{n+m+1}\s}^0=0$. By \cite[Lemma~3.3]{SZ1}, there is a uniform $B>0$ so that $\Re(g_{\kappa,\s'}^{m'}(t'))>t'-B$ for all $m'$, for all $\kappa$ from a bounded domain and all sufficiently large $t'$ (depending only on the bound $\kappa$), and for all $\s'$. Choosing $n$ sufficiently large, we can be sure that  $\Re\left( g_{\kappa,\sigma^{n+1}(\s)}^m(F^{\circ(n+1)}(t))-\kappa\right)>C$ for any given $C>2$ and all $m$, and this proves inductively that 
\[
\left|\frac{\partial}{\partial\kappa} g_{\kappa,\sigma^n(\s)}^{m+1}(F^{\circ n}(t))\right|
< \frac{2/C+1}{C} <2/C
\]
for all $m$ and thus also in the limit $m\to\infty$. The limit for $n\to\infty$ is thus equal to $0$ as claimed.

The remaining claim is $(d/d\kappa)E^n(\kappa)\to\infty$ as $n\to\infty$. For this, we use that $\kappa$ is on a parameter ray, say $\kappa=\Ray(t)$, and evaluate 
\[
dE^n(\kappa)/d\kappa=\lim_{\kappa_m\to\kappa}(E^n(\kappa_m)-E^n(\kappa))/(\kappa_m-\kappa) 
\]
using
$\kappa_m=\Ray(t_m)$ for a sequence $t_m\to t$. Since parameter rays $\kappa(t):=\Ray(t)$ are differentiable with $d\kappa/dt=\Ray'(t)\neq 0$, we can write 
\[
\frac{d}{d\kappa} E^n(\kappa) = \frac{d E^n(\Ray(t))}{dt} \cdot\frac{1}{\Ray'(t)}\,\,.
\]
Since we proved above that $(d/dt)E^n(\Ray(t))\to\infty$ as $n\to\infty$, this proves the last claim.
\qed

\proofof{Theorem~\ref{Thm:ParaRaysEndpoints}}
We only need to prove that for every open and bounded $\Lambda\subset\C$, we have $I_R\cap\Lambda\subset I_{p,\Lambda}$ for every $p>1$: once we know that, it follows $\dim_H(I_R\cap\Lambda)\le 1+1/p$ by Theorem~\ref{Thm:HausdorffDimParaRays}; by countable additivity and because this holds for all $p>1$, we have $\dim_H(I_R)\le 1$. But since $I_R$ contains curves, we conclude $\dim_H(I_R)=1$. Since $I=I_R\cup I_E$ and $\dim_H(I)=2$ by \cite{QiuDimTwo}, it follows that $\dim_H(I_E)=2$.

It remains to prove that $I_R\cap\Lambda\subset I_{p,\Lambda}$ for every $p>1$; more precisely, for $\kappa\in I_R$ we need to prove that $E^n(\kappa)\in P_{p,0}$ for sufficiently large $n$ and $|(E^n(\kappa))'|\to\infty$ as $n\to\infty$. The first statement is \cite[Proposition~4.5]{SZ1} and the second one is Lemma~\ref{Lem:ParaRaysDiffEstimate}.
\qed

\end{document}